\documentclass[a4paper,10pt]{article}
\usepackage{latexsym}
\usepackage{amsmath}
\usepackage{amscd}
\usepackage{amssymb}
\usepackage{url}
\usepackage{accents}

\usepackage[dvipdfm,%
 bookmarks=true,%
 bookmarksnumbered=true,%
 colorlinks=true,%
 pdftitle={self-emb},%
 pdfauthor={Yokoyama},%
 pdfsubject={},%
 pdfkeywords={}]{hyperref}

\usepackage{amsthm}

\def\WKL{\mathsf{WKL_0}}

\def\PA{\mathrm{PA}}
\def\CPA{\mathrm{CPA}}

\def\E{\exists}
\def\A{\forall}
\def\N{\mathbb{N}}
\def\Z{\mathbb{Z}}
\def\Q{\mathbb{Q}}

\def\rest{\upharpoonright}

\newcommand{\mr}[1]{\mathrm{#1}}
\newcommand{\mc}[1]{\mathcal{#1}}

\newcommand{\TP}[2]{\Sigma^{#1}_{#2}\text{-}\mathrm{TP}}

\def\P2{\Pi^1_2}

\def\TP1{\Pi^0_1\mathrm{TP}}

\newtheorem{thm}{Theorem}[section]
\newtheorem{theorem}[thm]{Theorem}

\newtheorem*{claim*}{Claim}
\newtheorem{prop}[thm]{Proposition}
\newtheorem{lem}[thm]{Lemma}

\newtheorem{corollary}[thm]{Corollary}

\theoremstyle{definition}
\newtheorem{defi}{Definition}[section]

\newtheorem{quest}[defi]{Question}
\newtheorem{exa}[defi]{Example}

\begin{document}

\title{A generalization of Levin-Schnorr's theorem}

\author{Keita Yokoyama%
\footnote{
School of Information Science, Japan Advanced Institute of Science and Technology,
1-1 Asahidai, Nomi, Ishikawa, 923-1292 Japan
e-mail: {\sf y-keita@jaist.ac.jp}
}
\footnote{
The author is grateful to Kenshi Miyabe for useful comments.
His work is partially supported by
 JSPS Grant-in-Aid for Research Activity Start-up grant number 25887026,
 JSPS-NUS Bilateral Joint Research Projects J150000618 (PI's: K.~Tanaka, C.~T.~Chong),
 JSPS fellowship for research abroad,
 and JSPS Core-to-Core Program (A.~Advanced Research Networks).
}
}
\date{}

\maketitle

\begin{abstract}
In this paper, we will generalize the definition of partially random or complex reals, and then show the duality of random and complex, \textit{i.e.}, a generalized version of Levin-Schnorr's theorem.
We also study randomness from the view point of arithmetic using the relativization to a complete $\Pi^{0}_{1}$-class.


\end{abstract}

\section{Introduction}
%

The notion of randomness is studied from several approaches.
Here, we would like to consider the major two approaches, namely, randomness defined by a measure, which is a generalization of Martin-L\"of randomness, and randomness defined by a complexity function, which is a generalization of weak Chaitin randomness.
(In this paper, we call the latter notion ``complex''.)
It is well-known, as Levin-Schnorr's theorem, that Martin-L\"of randomness and weak Chaitin randomness coincide.
(See, e.g., Downey and Hirschfeldt\cite{Downey-Hirschfeldt} or Nies\cite{Nies2009}.)
Then, is there a general correspondence between these two approaches?
In this paper, we will try to give a concrete connection between them.
We will provide some general definitions for the above two styles of randomness, and show a generalization of Levin-Schnorr's theorem. It will show that these two definitions of randomness have a duality, in other words, given a new notion of randomness in one of the above, then, one can automatically get the definition of the same notion of the other style.
In fact, our generalization of randomness cannot capture the whole known notions, but still cover several important notions of partial randomness (e.g. in Tadaki\cite{Tadaki} or Calude/Staiger/Terwijn\cite{Calude/Staiger/Terwijn}).
Based on our generalization, we also study on randomness in arithmetic using relativization to some Muchnik complete $\Pi^{0}_{1}$-class.

\section{Generalizing notions of complex and random}
In this section, we introduce a generalized notion of random or complex reals.
We first define random reals relative to a recursive sub-measure on all of codes for open sets.

\begin{defi}
 A \textit{pre-measure} is a recursive function $m:[2^{<\omega}]^{<\omega}\to [0,\infty)$ (more precisely, a recursive function $m:[2^{<\omega}]^{<\omega}\times \omega\to \Q$ such that $m(F,\cdot)$ codes a non-negative recursive real) which satisfies the following:
 
\begin{enumerate}
 \item $m(\emptyset)=0$,
 \item if $F_{1}\subseteq F_{2}$, then $m(F_{1})\le m(F_{2})$,
 \item $m(F_{1}\cup F_{2})\le m(F_{1})+m(F_{2})$.
\end{enumerate}
Given a pre-measure $m$ on all of codes for clopen sets, we expand it into a sub-measure on all of codes for open sets ${m}:[2^{<\omega}]^{\le \omega}\to [0,\infty)$ as follows: let $A\subseteq 2^{<\omega}$, then,
\[{m}(A)=\sup\{m(F)\mid F\subseteq_{\mr{fin}}A\}.\]

Let $Z\in 2^{\omega}$.
An \textit{${m}$-test (relative to $Z$)} is a uniformly ($Z$-)r.e.~sequence $\{A_{i}\mid i\in\omega\}$ such that $m(A_{i})\le2^{-i}$. A real $X\in\omega$ is said to be \textit{${m}$-random (relative to $Z$)} if $X\notin \bigcap_{i}[A_{i}]$ for any {${m}$-test (relative to $Z$)} $\{A_{i}\mid i\in\omega\}$.
(Here, $[A]$ is an open set generated by $A$, \textit{i.e.}, $[A]=\{X\in 2^{\omega}\mid \E \sigma\in A\  \E n\in\omega\  X\rest n=\sigma\}$.)
A \textit{Solovay-${m}$-test (relative to $Z$)} is a ($Z$-)r.e.~set $A\subseteq 2^{<\omega}$ such that $m(A)<\infty$. A real $X\in\omega$ is said to be \textit{Solovay-${m}$-random (relative to $Z$)} if for any {Solovay-${m}$-test (relative to $Z$)} $A$, there exists a finite subset $F\subseteq A$ such that $X\notin [A\setminus F]$.

\end{defi}
\begin{exa}
\begin{enumerate}
 \item Let $h:2^{<\omega}\to \omega$ be a recursive function. Then, the following are pre-measures:
 \begin{align*}
 \mr{dwt}_{h}(F)&:=\sum_{\sigma\in F}2^{-h(\sigma)},\\
 \mr{pwt}_{h}(F)&:=\sup\{\mr{dwt}_{h}(P)\mid \text{$P\subseteq_{\mr{fin}} F$ is prefix free}\},\\
 \mr{dct}_{h}(F)&:=\sup_{n\in\omega}\frac{\{\sigma\in F\mid h(\sigma)<n\}^{\#}}{2^{n}},\\
 \mr{pct}_{h}(F)&:=\sup\{\mr{dct}_{h}(P)\mid \text{$P\subseteq_{\mr{fin}} F$ is prefix free}\}.
\end{align*}
Here, $\mr{dwt}_{h}$-random is usually called $h$-random, which is appeared, e.g., in Tadaki\cite{Tadaki}, and, in particular, it is Martin-L\"of random if $h(\sigma)=|\sigma|$.
$\mr{pwt}_{h}$-random is usually called strongly-$h$-random, which is appeared, e.g., in Calude/Staiger/Terwijn\cite{Calude/Staiger/Terwijn}.
Note that the original notion of Martin-L\"of random and Solovay random are equivalent, but they are different in case, e.g., $m=\mr{dwt}_{h}$.
 \item If $m_{1}$ and $m_{2}$ are pre-measures, then, $m_{1}+m_{2}$ and $\max\{m_{1},m_{2}\}$ are pre-measures.
  \item In fact, any $\Sigma^{0}_{2}$-subclass of Cantor space can be considered as a set of $m$-random reals, and conversely, for any pre-measure $m$, a class of $m$-random real is a $\Pi^{0}_{2}$-class.
  Let $P\subseteq 2^{\omega}$ be a $\Sigma^{0}_{2}$-class. Take a recursive sequence of trees $\{ T_{i}\mid i\in\omega\}$ such that $X\in P$ if and only if $X$ is a path of $T_{i}$ for some $i\in \omega$.
  Define (recursive) pre-measures $m_{i}$ and $m$ as follows: 
\begin{align*}
 m_{i}(F)&=
\begin{cases}
 1 & \mbox{if $F\cap T_{i}\neq\emptyset$},\\
 0 & \mbox{otherwise},
\end{cases}
\\
m(F)&=\sum_{i\in\omega}2^{-i}m_{i}(F).
\end{align*}
Then, we can easily check that $X$ is $m$-random if and only if $X\in P$.
\end{enumerate}
\end{exa}
\begin{prop}
 For any pre-measure $m$ and for any $Z\in 2^{\omega}$, a universal ${m}$-test exists, in other words, there exists an {${m}$-test relative to $Z$} $\{A_{i}\mid i\in\omega\}$ such that $X\in2^{\omega}$ is $m$-random relative to $Z$ if and only if $X\notin\bigcap_{i}[A_{i}]$.
\end{prop}
\begin{proof}
Similar to the construction of a universal Martin-L\"of test.
\end{proof}

Next, we define generalized complexities in two ways.
We generalize two different style definitions introduced by Uspensky and Shen \cite{Uspensky-Shen}.
We first generalize the complexity defined by a description mode.
This provides a natural generalization of $\mr{KP}$, $\mr{KS}$, $\mr{KM}$ or $\mr{KD}$.
\begin{defi}[Complexity defined by a description mode]
 A \textit{rule} for a description mode is a recursive set $R\subset [2^{<\omega}\times2^{<\omega}]^{<\omega}$ which satisfies the following:
\begin{enumerate}
 \item $\emptyset\in R$.
 \item If $r\in R$ and $s\subseteq r$, then $s\in R$.
 \item If $r,s\in R$, then $\{(0^{\frown}\tau,\sigma)\mid (\tau,\sigma)\in r\}\cup\{(1^{\frown}\tau,\sigma)\mid (\tau,\sigma)\in s\}\in R$.
\end{enumerate}
Let $Z\in2^{\omega}$.
 A \textit{mode (relative to $Z$)} is a ($Z$-)r.e.~set $M\subseteq 2^{<\omega}\times 2^{<\omega}$, and we define the $M$-complexity $K^{M}:2^{<\omega}\to\N$ as $K^{M}(\sigma)=\min\{|\tau|\mid (\tau,\sigma)\in M\}$.
A mode $M$ is said to be an $R$-mode if any finite subset of $M$ is a member of $R$, and an $R$-mode $M$ is said to be \textit{$R$-optimal} if for any $R$-mode $M'$, there exists $c^{M'}\in\omega$ such that $K^{M}(\sigma)\le K^{M'}(\sigma)+c^{M'}$.
If $M$ is $R$-optimal, $K^{M}$ is called $R$-complexity.

Let $M$ be an $R$-optimal mode (relative to $Z$), then $X\in 2^{\omega}$ is said to be \textit{$R$-complex (relative to $Z$)} if there exists $c\in\omega$ such that for any $n\in\omega$, $K^{M}(X\rest n)\ge n-c$.
\end{defi}
\begin{exa}
 Prefix-free complexity $\mr{KP}$, simple complexity $\mr{KS}$, monotone complexity $\mr{KM}$ and decision complexity $\mr{KD}$ can be defined in this way.
For example, $R_{\mr{KP}}=\{r\in [2^{<\omega}\times2^{<\omega}]^{<\omega}\mid r$ is a finite partial function from $2^{<\omega}$ to $2^{<\omega}$ whose domain is prefix-free$\}$ is the rule for $\mr{KP}$.
See Uspensky/Shen \cite{Uspensky-Shen}. (In their paper, they use the word `entropy' in stead of `complexity'.)
\end{exa}
\begin{prop}
For any rule $R$ and for any $Z\in2^{\omega}$, an $R$-optimal mode relative to $Z$ exists.
\end{prop}
\begin{proof}
Similar to the construction of optimal prefix-free Turing machine. 
\end{proof}

Next, we introduce the definition of complexity as a minimal function again following the idea by Uspensky and Shen \cite{Uspensky-Shen}.
Here, we will through away some information which can be captured by a mode, e.g., we cannot describe $\mr{KM}$ in this way.
The definition is rather tricky, but in fact, this is the right notion corresponding to the randomness defined by sub-measure in the sense of Levin-Schnorr's theorem.

A finite complexity function is a finite set $r\subseteq 2^{<\omega}\times\Z$, we identify $r$ as a function $K^{r}(\sigma)=\min\{d\mid (\sigma,d)\in r\}\cup\{\infty\}$. Given a finite complexity $r\subseteq 2^{<\omega}\times\Z$, define $\ring{r}=\{\sigma\in2^{<\omega}\mid \E d\in\omega\  (\sigma,d)\in r\}$ and $r^{+i}:=\{(\sigma,d+i)\mid (\sigma,d)\in r\}$.
Let $r,s\subseteq 2^{<\omega}\times\omega$ be finite complexity functions, we say that $r$ is stronger than $s$ ($s\prec r$) if for any $(\sigma,d)\in s$, there exists $d'\le d$ such that $(\sigma,d')\in r$.
\begin{defi}[Complexity as a minimal function]
 A \textit{rule} (for a complexity function) is a recursive set $R\subset [2^{<\omega}\times\Z]^{<\omega}$ which satisfies the following:
\begin{enumerate}
 \item $\emptyset\in R$.
 \item If $r\in R$ and $s\prec r$, then $s\in R$.
 \item If $r,s\in R$, then $(r\cup s)^{+1}\in R$.
\end{enumerate}
 A \textit{complexity function (relative to $Z$)} is a right ($Z$-)r.e. function $K:2^{<\omega}\to \omega$.
(Here, we say that $K$ is a right r.e.~function if the relation $\{(\sigma,m)\mid K(\sigma)<m\}$ is r.e.)
 Given a rule $R$, a complexity function $K=K_{R}$ ($K=K^{Z}_{R}$) is said to be \textit{$R$-optimal (relative to $Z$)} if
\begin{enumerate}
 \item $R$-function: for any finite $F\subseteq 2^{<\omega}$, $\{(\sigma,K(\sigma))\mid \sigma\in F\}\in R$.
 \item $R$-minimal: if $A\subseteq 2^{<\omega}\times \omega$ is a ($Z$-)r.e. set such that any finite $F\subseteq A$ is an element of $R$, then there exists $c\in\omega$ such that for any $(\sigma,d)\in A$, $K(\sigma)<d+c$. 
\end{enumerate}
$X\in2^{\omega}$ is said to be $R$-complex relative to $Z$ if there exists $c\in\omega$ such that $K^{Z}_{R}(X\rest n)\ge n-c$ for any $n\in\omega$.
$X\in2^{\omega}$ is said to be Solovay-${R}$-complex relative to $Z$ if $\lim_{n\to\infty}K^{Z}_{R}(X\rest n)-n=\infty$.
Note that, sometimes, $R$-optimal complexity function is just called $R$-complexity.
\end{defi}

\begin{prop}
For any rule $R$ and for any $Z\in2^{\omega}$, $R$-optimal complexity function $K^{Z}_{R}$ exists.
\end{prop}
\begin{proof}
Let $\{A_{i}\mid i\in\omega\}$ be a recursive enumeration of all ($Z$-)r.e.~sets such that any finite $F\subseteq A_{i}$ is an element of $R$.
Define $\tilde{A}=\bigcup_{i}A_{i}^{+i}$, and define $K^{Z}_{R}(\sigma)=\min\{d\mid (\sigma,d)\in \tilde{A}\}\cup\{\infty\}$.
Then, we can easily check that this $K^{Z}_{R}$ is a desired function.
\end{proof}

\begin{exa}
\begin{enumerate}
 \item Let $h:2^{<\omega}\to \omega$ be a recursive function.
Then, the following are rules:
\begin{align*}
R_{\mr{KP}{h}}=&\left\{r\mid \sum_{ (\sigma,d)\in r}2^{-d+|\sigma|-h(\sigma)}<1\right\},\\
R_{\mr{KA}{h}}=&\left\{r\mid \sum_{ (\sigma,d)\in s}2^{-d+|\sigma|-h(\sigma)}<1\text{ for any $s\subseteq r$ such that $\ring{s}$ is prefix-free}\right\},\\
R_{\mr{KS}{h}}=&\left\{r\mid \{(\sigma,d)\in r\mid d-|\sigma|+h(\sigma)<n\}^{\#}<2^{n}\text{ for any $n\in\omega$}\right\},\\
R_{\mr{KD}{h}}=&\left\{r\mid \{(\sigma,d)\in s\mid d-|\sigma|+h(\sigma)<n\}^{\#}<2^{n}\right.\\
&\quad \text{ for any $n\in\omega$ and for any $s\subseteq r$ such that $\ring{s}$ is prefix-free}\}.
\end{align*}
Then, $K_{R_{\mr{KP}{h}}}(\sigma)=KP(\sigma)-|\sigma|+h(\sigma)$, $K_{R_{\mr{KA}{h}}}(\sigma)=KA(\sigma)-|\sigma|+h(\sigma)$, $K_{R_{\mr{KS}{h}}}(\sigma)=KS(\sigma)-|\sigma|+h(\sigma)$ and $K_{R_{\mr{KD}{h}}}(\sigma)=KD(\sigma)-|\sigma|+h(\sigma)$ up to constant, respectively, where $\mr{KA}$ is a priori complexity.
$R_{\mr{KP}{h}}$-complex is usually called $h$-complex, and, in particular, it is called weak Chaitin random if $h(\sigma)=|\sigma|$.
$R_{\mr{KA}{h}}$-complex is usually called strongly-$h$-complex.
\item If $R_{1}$ and $R_{2}$ are rules, then $R_{1}\cap R_{2}$ and $R_{1}\cup R_{2}\cup \{(r\cup s)^{+1}\mid r\in R_{1},s\in R_{2}\}$ are rules.
\end{enumerate}

\end{exa}

If $R$ is a rule for a mode, and $r\in R$, define $\hat{r}=\{(\sigma,|\tau|)\mid (\tau,\sigma)\in r\}$, then $\hat{R}:=\{s\mid \E r\in R\  s\prec \hat{r}\}$ is a rule for a complexity function.
If $M$ is an $R$-mode, then $K^{M}$ is an $\hat{R}$-complexity function.



\section{Generalized Levin-Schnorr's theorem}
In this section, we will show that randomness defined by a measure and complex defined by a complexity function have a concrete correspondence.
In this section, a rule means a rule for a complexity function $R\subset [2^{<\omega}\times\Z]^{<\omega}$.
For $r\in [2^{<\omega}\times\Z]^{<\omega}$, define $\|r\|=\min\{|\sigma|-d\mid (\sigma,d)\in r\}\in \Z\cup\{\infty\}$ ($\|\emptyset\|=\infty$).
We can easily check that $s\prec r$ if $\ring{s}\subseteq\ring{r}$ and $\|s\|\le \|r\|$.
\begin{lem}\label{union-of-rules}
 Let $R$ be a rule, and let $r_{1},\dots,r_{n}\in R$.
 Then, $r_{1}^{+1}\cup\dots\cup r_{n}^{+n}\in R$.
\end{lem}
\begin{proof}
 By induction on $n$.
 If $r_{2}^{+1}\cup\dots\cup r_{n}^{+n-1}\in R$ and $r_{1}\in R$, then, $(r_{1}\cup(r_{2}^{+1}\cup\dots\cup r_{n}^{+n-1}))^{+1}=r_{1}^{+1}\cup\dots\cup r_{n}^{+n}\in R$.
\end{proof}
\begin{defi}
 Let $m$ be a pre-measure, and let $R$ be a rule.
 Then, we define $m^{\surd}\subseteq [2^{<\omega}\times\omega]^{<\omega}$ and $R^{\surd}:[2^{<\omega}]^{<\omega}\to [0,\infty)$ as follows:
\begin{align*}
 m^{\surd}&:=\{r\in [2^{<\omega}\times\omega]^{<\omega}\mid \A s\subseteq r\  m(\ring{s})\le2^{-\|s\|}\},\\
 R^{\surd}(F)&:=\inf\{2^{-\|r_{1}\|}+\dots+2^{-\|r_{l}\|}\mid r_{1},\dots,r_{l}\in R,\  F\subseteq \ring{r_{1}}\cup\dots\cup\ring{r_{l}}\}.
\end{align*}
\end{defi}
Note that in the above definition, the rule $m^{\surd}$ is essentially defined by the logarithm of measure $m$ as for the usual correspondence of measure and complexity.
\begin{prop}
 If $R$ is a rule, then $R^{\surd}$ is a pre-measure.
 If $m$ is a pre-measure, then $m^{\surd}$ is a rule.
\end{prop}
\begin{proof}
Let $R$ be a rule.
Then, $R^{\surd}(\emptyset)=2^{-\infty}=0$, monotonicity of $R^{\surd}$ is obvious from the definition, and
\begin{align*}
R^{\surd}(F_{1}\cup F_{2})&=\inf\{2^{-\|r_{1}\|}+\dots+2^{-\|r_{l}\|}\mid r_{1},\dots,r_{l}\in R,\  F_{1}\cup F_{2}\subseteq \ring{r_{1}}\cup\dots\cup\ring{r_{l}}\}\\
&\le \inf\{2^{-\|r_{1}\|}+\dots+2^{-\|r_{l}\|}\mid r_{1},\dots,r_{l}\in R,\  F_{1}\subseteq \ring{r_{1}}\cup\dots\cup\ring{r_{l}}\}\\
&\quad\quad +\inf\{2^{-\|r_{1}\|}+\dots+2^{-\|r_{l'}\|}\mid r_{1},\dots,r_{l'}\in R,\  F_{2}\subseteq \ring{r_{1}}\cup\dots\cup\ring{r_{l'}}\}\\
&=R^{\surd}(F_{1})+R^{\surd}(F_{2}).
 \end{align*}
Thus, $R^{\surd}$ is a pre-measure.

Let $m$ be a pre-measure. Then, by definition, $\emptyset\in m^{\surd}$ and $s\prec r\wedge r\in m^{\surd}$ implies $s\in m^{\surd}$.
Let $r_{1},r_{2}\in m^{\surd}$ and $s\subseteq (r_{1}\cup r_{2})^{+1}$.
Define $s_{1}=s^{-1}\cap r_{1}$ and $s_{2}=s^{-1}\cap r_{2}$, then, $\ring{s}=\ring{s_{1}}\cup\ring{s_{2}}$ and $\|s\|=\|(s_{1}\cup s_{2})^{+1}\|=\|s_{1}\cup s_{2}\|-1=\min\{\|s_{1}\|,\|s_{2}\|\}-1$.
Hence,
\begin{align*}
 m(\ring{s})=m(\ring{s_{1}}\cup\ring{s_{2}})\le m(\ring{s_{1}})+m(\ring{s_{2}})\le 2^{-\|s_{1}\|}+2^{-\|s_{2}\|}\le 2\cdot 2^{-\|s_{1}\cup s_{2}\|}=2^{-\|s\|}.
\end{align*}
Thus, $(r_{1}\cup r_{2})^{+1}\in m^{\surd}$.
We have proved that $m^{\surd}$ is a rule.
\end{proof}
\begin{prop}Let $m,k$ be pre-measures, and $R,S$ be rules.
\begin{enumerate}
\item If $m\le ck$ for some $c\in\omega$, then there exists $c'\in\omega$ such that $K_{m^{\surd}}\le K_{k^{\surd}}+c'$.
\item If $K_{R}\le K_{S}+c$ for some $c\in\omega$, then there exists $c'\in\omega$ such that $R^{\surd}\le c'S^{\surd}$.
\end{enumerate}
\end{prop}
\begin{proof}
 Easy from the definition.
\end{proof}
The following proposition means that $m$ and $m^{\surd\surd}$ is essentially the same, and $R$ and $R^{\surd\surd}$ is essentially the same.
\begin{prop}\label{prop-duality}
 Let $m$ be a pre-measure, and let $R$ be a rule.
Then,
\begin{enumerate}
\item  $m\le m^{\surd\surd}\le 2m$,
\item  $R\subseteq R^{\surd\surd}\subseteq\{s\mid \E r\in R\  s\prec r^{-2}\}$, thus, $K_{R}-c \le K_{R^{\surd\surd}}\le K_{R}+c$ for some $c\in\omega$.
\end{enumerate}
\end{prop}
\begin{proof}
We first prove 1.
 Let $m$ be a pre-measure, and $F$ be a finite subset of $2^{<\omega}$. Then,
\begin{align*}
 m^{\surd\surd}(F)&=\inf\{2^{-\|r_{1}\|}+\dots+2^{-\|r_{l}\|}\mid r_{1},\dots,r_{l}\in m^{\surd},\  F\subseteq \ring{r_{1}}\cup\dots\cup\ring{r_{l}}\}\\
&=\inf\{2^{-\|r_{1}\|}+\dots+2^{-\|r_{l}\|}\mid \A s\subseteq r_{i}(m(\ring{s})\le 2^{-\|s\|}),\  F\subseteq \ring{r_{1}}\cup\dots\cup\ring{r_{l}}\}\\
&\ge\inf\{m(\ring{r_{1}})+\dots+m(\ring{r_{l}})\mid F\subseteq \ring{r_{1}}\cup\dots\cup\ring{r_{l}}\}\\
&\ge m(F).
\end{align*}
We next show that $m^{\surd\surd}(F)\le 2m(F)$.
For $e\in \Z$, define $r_{e,F}=\{(\sigma,|\sigma|-e)\mid \sigma\in F\}$.
Then, for any $e\in \Z$ such that $m(F)\le 2^{-e}$ and for any non-empty $s\subseteq r_{e,F}$, we have $m(\ring{s})\le m(F)\le 2^{-e}= 2^{-\|s\|}$, thus, $r_{e,F}\in m^{\surd}$.
Hence,
\begin{align*}
 m^{\surd\surd}(F)&\le \inf\{2^{-\|r\|}\mid r\in m^{\surd},\  F\subseteq \ring{r}\}\\
&\le \inf\{2^{-\|r_{e,F}\|}\mid e\in\Z, m(F)\le 2^{-e}\}\\
&\le \inf\{2^{-e}\mid e\in Z, m(F)\le 2^{-e}\}\le 2m(F).
\end{align*}

Next, we prove 2.
Let $R$ be a rule for a complexity. Then,
\[r\in R^{\surd\surd}\leftrightarrow \A s\subseteq r\  \inf\{2^{-\|t_{1}\|}+\dots+2^{-\|t_{l}\|}\mid t_{1},\dots,t_{l}\in R,\  \ring{s}\subseteq \ring{t_{1}}\cup\dots\cup\ring{t_{l}}\}\le 2^{-\|s\|}.\]
Thus, $r\in R$ implies $r\in R^{\surd\surd}$.
Let $r\in R^{\surd\surd}$. Then, there exist $t_{1},\dots,t_{l}\in R$ such that $\ring{r}\subseteq \ring{t_{1}}\cup\dots\cup\ring{t_{l}}$ and $2^{-\|t_{1}\|}+\dots+2^{-\|t_{l}\|}\le 2^{-\|r\|+1}$. If $\|t_{i}\|=\|t_{j}\|$, we can choose $(t_{i}\cup t_{j})^{+1}$ in stead of $t_{i}$ and $t_{j}$ since $2^{-\|t_{i}\|}+2^{-\|t_{j}\|}=2^{-\|t_{i}\cup t_{j}\|+1}=2^{-\|(t_{i}\cup t_{j})^{+1}\|}$. Thus, without loss of generality, we can assume that $\|t_{1}\|<\|t_{2}\|<\dots<\|t_{l}\|$.
By Lemma~\ref{union-of-rules}, $t:=t_{1}^{+1}\cup\dots\cup t_{l}^{+l}\in R$.
Then, $\|t\|=\|t_{1}\|-1$ and $\ring{r}\subseteq\ring{t}$.
By $2^{-\|t_{1}\|}\le 2^{-\|r\|+1}$, we have $\|t\|\ge \|r\|-2$.
Thus, $r\prec t^{-2}$.
\end{proof}
\begin{defi}
 Let $m$ be a pre-measure, and let $R$ be a rule.
Then, $R$ is said to be a dual rule of $m$, or $m$ is said to be a dual pre-measure of $R$, if there exists $c\in\omega$ such that $K_{R}-c\le K_{m^{\surd}}\le K_{R}+c$, or equivalently by Proposition~\ref{prop-duality}, there exists $c\in\omega$ such that $1/c\cdot m\le R^{\surd}\le cm$. 
\end{defi}
\begin{exa}Let $h:2^{<\omega}\to \omega$ be a recursive function. Then,
\begin{itemize}
\item $R_{\mr{KP}{h}}$ is a dual of $\mr{dwt}_{h}$,
\item $R_{\mr{KA}{h}}$ is a dual of $\mr{pwt}_{h}$,
\item $R_{\mr{KS}{h}}$ is a dual of $\mr{dct}_{h}$,
\item $R_{\mr{KD}{h}}$ is a dual of $\mr{pct}_{h}$.
\end{itemize}
\end{exa}
\begin{theorem}[Duality/Generalized Levin-Schnorr's theorem]\label{general-Schnorr}
 Let $m$ be a pre-measure, and let $R$ be its dual rule.
Then, $X\in 2^{\omega}$ is $m$-random if and only if it is $R$-complex.
\end{theorem}
\begin{proof}
Without loss of generality, we may assume that $m=R^{\surd}$.
 Let $X\in 2^{\omega}$ be not $R$-complex.
 Define $U_{i}=\{\sigma\in 2^{<\omega}\mid K_{R}(\sigma)\le |\sigma|-i\}$.
 Then, a sequence $\{U_{i}\mid i\in\omega\}$ is (uniformly) r.e., and $X\in \bigcap_{i\in\omega} [U_{i}]$.
We show that ${m}(U_{i})\le 2^{-i}$.
Let $F\subseteq_{\mr{fin}}U_{i}$. Define $r_{F}=\{(\sigma,K_{R}(\sigma))\mid \sigma\in F\}$.
Then, $r_{F}\in R$ and $\|r\|\ge i$. Thus, $m(F)=R^{\surd}(F)\le 2^{-\|r_{F}\|}\le 2^{-i}$.
Hence, $X$ is not $m$-random.

Conversely, let $X\in 2^{\omega}$ be not $m$-random.
Then, there exists an ${m}$-test $\{U_{i}\mid i\in\omega\}$ such that $X\in \bigcap_{i\in\omega}[U_{i}].$
Define an r.e.~set $A$ as $A=\{(\sigma,|\sigma|-i)\mid \sigma\in U_{2i}\}$.
Let $a\subseteq_{\mr{fin}}A$, and $s\subseteq a$.
If $\|s\|=n$, then $\ring{s}\subseteq \bigcup_{i\ge n}U_{2i}$, thus,
\[ m(\ring{s})\le {m}\left(\bigcup_{i\ge n}U_{2i}\right)\le 2^{-n}=2^{-\|s\|}.\]
Hence, $a\in m^{\surd}=R^{\surd\surd}$.
Thus, $X$ is not $R^{\surd\surd}$-complex, and hence it is not $R$-complex by Proposition~\ref{prop-duality}.
\end{proof}
\begin{theorem}[Duality/Generalized Levin-Schnorr's theorem]
 Let $m$ be a pre-measure, and let $R$ be its dual rule.
Then, $X\in 2^{\omega}$ is Solovay-$m$-random if and only if it is Solovay-$R$-complex.
\end{theorem}
\begin{proof}
 Similar to Theorem~\ref{general-Schnorr}.
\end{proof}

\section{Random and complex relative to PA degree \\ and nonstandard models of arithmetic}
In this section we will show that a real $X$ is complex/random from a nonstandard model of arithmetic if and only if it is not compressible in arithmetic, if and only if it is `strongly' complex/random.
An example of this is that $X$ is Martin-L\"of random if and only if it is Martin-L\"of random from a nonstandard model of arithmetic if and only if it is not compressible in a prefix-free way in arithmetic.
(This is an easy consequence of Theorem~\ref{gen-random-relPA} or Corollary~\ref{cor-complex-in-arith}.)
We will see the meaning of this in a general setting based on the previous section.
For this,  we consider the concept complex/random relative to some Muchnik complete $\Pi^{0}_{1}$-class.

In this section, we fix a Muchnik complete $\Pi^{0}_{1}$-class $\CPA\subseteq 2^{\omega}$.
(Here, $\CPA$ means a class of completions of Peano Arithmetic, which is a well-known Muchnik/Medvedev complete $\Pi^{0}_{1}$-class.)
We say that $X\in 2^{\omega}$ is $m$-random (or $R$-complex) relative to $\CPA$ if there exists $Z\in\CPA$ such that $X$ is $m$-random (or $R$-complex) relative to $Z$.
We fix a theory of arithmetic $T=\PA$ or a recursive extension of  $\mr{I}\Sigma_{1}$.
Let $R$ be a $\Sigma_{0}$-definable rule (or $T$ provably $\Delta_{1}$ rule).

\begin{defi}
 Let $X\in 2^{\omega}$, and let $M$ a model of $T$.
 For a given $r\in M$ such that $M\models r\in R$, define $K_{r}:2^{<M}\to M\cup\{\infty\}$ as $K_{r}(\sigma)=\min\{d\mid (\sigma,d)\in r\}\cup\{\infty\}$ ($r$ or $K_{r}$ is said to be an $M$-finite complexity).
Then, $X$ is said to be $R^{M}$-complex if for any $M$-finite complexity $r\in R$ in $M$, there exists $c\in\omega$ such that $M\models K_{r}(X\rest n)\ge n-c$ for any $n\in\omega$.
\end{defi}

\begin{theorem}\label{complex-nonst}
 Let $X\in 2^{\omega}$. Then, the following are equivalent.
\begin{enumerate}
 \item $X$ is $R$-complex relative to $\CPA$.
 \item $X$ is $R^{M}$-complex for some nonstandard model $M\models T$.
\end{enumerate}
\end{theorem}
\begin{proof}
 We first show 1 $\to$ 2.
 Let $X$ be $R$-complex relative to $\CPA$.
 Then, there exists a Scott set $S\subseteq \mc{P}(\omega)$ (\textit{i.e.}, $(\omega, S)\models \WKL$) such that for any $A\in S$, $X$ is $R$-complex relative to $A$.
 By Theorem~15.23 of \cite{Kaye}, there exists a nonstandard model $M\models T$ such that $\mr{SSy}(M)=S$.
 Then, we can easily check that $X$ is $R^{M}$-complex.
 
 To show 2 $\to$ 1, we fix a $\Sigma_{0}$-formula $\theta(\sigma,n,m,\tau)$ such that for any $A\in 2^{\omega}$, $K^{A}(\sigma)=\min\{n\mid \E m\in\omega\  \omega\models\theta(\sigma,n,m,A\rest m)\}$ be an optimal $R$-complexity relative to $A$, and $T$ proves `for any finite $r$ such that $(\sigma,n)\in r\to \E m\theta(\sigma,n,m,\rho\rest m)$ for some $\rho$, $r$ is a finite $R$-complexity'.
 Let $M\models T$ and $X$ be $R^{M}$-complex.
 Then, there exists $A\in\CPA$ such that $A\in\mr{SSy}(M)$. Take $\rho\in M$ such that $\mr{SSy}(\rho)=A$, and, in $M$ define $K^{\rho}:2^{<M}\to M$ as $K^{\rho}(\sigma)=\min\{n\mid \E m\le|\rho|\theta(\sigma,n,m,\rho\rest m)\}\cup\{\infty\}$.
 Then, for some $H\in M\setminus\omega$, $K^{\rho}\rest H$ is a finite $R$-complexity in $M$. Thus, there exists $c\in\omega$ such that for any $n\in\omega$, $K^{A}(X\rest n)\ge K^{\rho}(X\rest n)\ge n-c$.
 Hence, $X$ is $R$-complex relative to $A$.
\end{proof}

Theorem~\ref{complex-nonst} shows that if $X\in2^{\omega}$ is an $R$-complex relative to $\CPA$, then it is not $R$-compressible in $T$ in the following sense. 
Let $K$ be a new function symbol. Define $T^{R}:=\PA(K)+$`$K$ is an $R$-complexity function' (in other words, any finite part of $K$ is in $R$), and define ${T^{R*}}=T^{R}+\{K(\sigma)=\min\{n\mid \E m\in\omega\  \omega\models\theta(\sigma,n,m,0^{m})\}\}$ where $\theta$ is a $\Sigma_{0}$-formula defining $R$-optimal complexity appeared in the proof of Theorem~\ref{complex-nonst}.
For $T'\supseteq T^{R}$ and $X\in 2^{\omega}$, $X$ is said to be \textit{compressible in $T'$} if for any $c\in\omega$, there exists $n\in\omega$ such that $T'\vdash K(X\rest n)<n-c$.
Then, the following is an easy modification of Theorem~\ref{complex-nonst}.

Let $\pi(e,X,\sigma,n)\equiv \E m \pi_{0}(m,e,X[m],\sigma,n)$ be a $\Sigma^{0}_{1}$-universal lightface formula, \textit{i.e.}, for any $\Sigma^{0}_{1}$-formula $\varphi(X, \sigma,n)$, there exists $e_{0}\in \omega$ such that $\mr{I}\Sigma^{0}_{1}$ proves $\pi(e_{0},X,\sigma,n)\leftrightarrow\varphi(X, \sigma,n)$.
Let $S^{X,e}:=\{(\sigma,n)\mid \E N((\sigma,n)\in \{(\sigma',n')<N\mid \E m<N\pi_{0}(e,X[m],\sigma',n')\}\in R)\}$, and define $K^{X}_{R}(\sigma)=\min\{n-e\mid (\sigma,n)\in S^{X,e}\}$.
This $K^{X}_{R}$ is a $\Pi^{X}_{2}$-definable (actually, $\Sigma^{X}_{1}\wedge \Pi^{X}_{1}$-definable) function.
For $A\in 2^{\omega}$, we can easily check that $K^{A}_{R}$ is an $R$-optimal complexity function relative to $A$.
In general, we can show the following.
Let $(M,A)\models \mr{I}\Sigma^{0}_{1}$, and let $\varphi(\sigma,n,A)$ be a $\Sigma_{1}^{A}$-formula without parameters such that $(M,A)\models \A N \{(\sigma,n)<N\mid \varphi(\sigma,n,A) \}\in R)$.
Then, there exists $c\in\omega$ such that $(M,A)\models\A \sigma \A n(\varphi(\sigma,n,A)\to K^{A}_{R}(\sigma)<n+c)$.

\begin{lem}
 Let $\{\sigma_{i}\}_{i\in \omega}\subseteq 2^{<\omega}$, and let $\{n_{i}\}_{i\in \omega}\subseteq \omega$.
 Then, the following are equivalent.
 
\begin{enumerate}
 \item There exists an consistent recursive extension $T'\supseteq T$ such that for any $c\in \omega$ there exists $l\in \omega$ such that $T'\vdash \bigvee_{i<l}(K_{R}(\sigma_{i})<n_{i}-c)$.
 \item For any $A\in \CPA$, and for any $c\in\omega$, there exists $i\in\omega$ such that $K_{R}^{A}(\sigma_{i})<n_{i}-c$.
\end{enumerate}
\end{lem}
\begin{proof}
 To show 1 $\to$ 2, let $T'$ be such that for any $c\in \omega$ there exists $l\in \omega$ such that $T'\vdash \bigvee_{i<l}(K_{R}(\sigma_{i})<n_{i}-c)$, and let $A\in\CPA$.
 Then, as in the proof of Theorem~\ref{complex-nonst}, there exists a countable model $M\models T'$ such that for any $X\in \mr{SSy}(M)$, $X$ is Turing reducible to $A$.
 Thus, the set $S:=\{(\sigma,n)\in 2^{<\omega}\times\omega \mid M\models K_{R}(\sigma)\le n\}$ is Turing reducible to $A$, and any finite subset of $S$ is a member of $R$.
 By the optimality of $K_{R}^{A}$, there exists $c_{0}\in\omega$ such that for any $(\sigma,n)\in S$, $K_{A}^{R}(\sigma)<n+c_{0}$.
Let $c\in\omega$.
Then, by the assumption of $T'$, there exists $i\in\omega$ such that $(\sigma_{i},n_{i}-c-c_{0})\in S$ for any $i\in\omega$.
Thus, $K_{R}^{A}(\sigma_{i})<n_{i}-c$.
 
Next, we will show $\neg$1 $\to$ $\neg$2.
Let $T'=T+\mr{Con}(T)$.
Then, there exists $c\in \omega$ such that $T'\not\vdash\bigvee_{i<l}(K_{R}(\sigma_{i})<n_{i}-c)$ for any $l\in\omega$, thus, $T'+\{K_{R}(\sigma_{i})\ge n_{i}-c\mid i\in\omega\}$ is consistent.
Fix a $\Sigma_{0}$-definable predicate $\mc{X}\subseteq 2^{<\omega}$ such that the class of paths of $\mc{X}$ is a Muchnik complete $\Pi^{0}_{1}$-class, and $T$ proves $\A n\E \sigma(|\sigma|=n\wedge \sigma\in \mc{T})$ (e.g., take $\mc{X}$ such that its path is a completion of $\mr{I}\Sigma_{0}$).
We will show that for some path $A\in 2^{\omega}$ of $\mc{X}$ and for some $C\in\omega$, $K^{A}_{R}(\sigma_{i})\ge n_{i}-C$ for any $i\in\omega$.
Let $M\models T'+\{K_{R}(\sigma_{i})\ge n_{i}-c\mid i\in\omega\}$, let $\lambda=\min\{m\in M\mid$ $m$ is a proof of $T\vdash\bot\}$, and let $\tau$ be the leftmost path branch of $\mc{X}^{=\lambda}$ in $M$.
Note that $\tau$ is a $\Sigma_{0}$-definable element in $M$, and $|\tau|>m$ for any $m\in\omega$. 
Within $M$, a set $S:=\{(\sigma,m)\mid K_{R}^{\tau}(\sigma)\le m\}$ is $\Sigma_{1}$-definable without any parameters. Thus, there exists $c_{0}\in\omega$ such that $M\models\A(\sigma,m)\in S(K_{R}(\sigma)<m+c_{0})$. Define $A=\tau\rest \omega\in 2^{\omega}$. Then, $A\in \mc{X}$ (in the standard model).
Let $d_{i}=K_{R}^{A}(\sigma_{i})$ (in the standard model). Then, within $M$, $d_{i}\ge K_{R}^{\tau}(\sigma_{i})$, \textit{i.e.}, $(\sigma_{i},d_{i})\in S$.
Therefore, we have $M\models d_{i}+c_{0}>K_{R}(\sigma_{i})\ge n_{i}-c$.
Hence, $d_{i}=K^{A}_{R}(\sigma_{i})\ge n_{i}-c-c_{0}$ for any $i\in\omega$.
\end{proof}
\begin{lem}
 Let $\{\sigma_{i}\}_{i\in \omega}\subseteq 2^{<\omega}$, and let $\{n_{i}\}_{i\in \omega}\subseteq \omega$, and let $k\in\omega$.
 Then, the following are equivalent.
 
\begin{enumerate}
 \item There exists an consistent recursive extension $T'\supseteq T$ such that for any $c\in \omega$ there exists $l\in \omega$ such that $T'\vdash \bigvee_{i<l}(K^{\mathbf{0}^{(k)}}_{R}(\sigma_{i})<n_{i}-c)$.
 \item For any $A\in \CPA$, and for any $c\in\omega$, there exists $i\in\omega$ such that $K_{R}^{A}(\sigma_{i})<n_{i}-c$.
\end{enumerate}
\end{lem}
\begin{proof}
 Similar to the previous lemma.
 (This time, $K^{\mathbf{0}^{(k)}}_{R}$ is a $\Pi_{k+2}$-definable function. The set $S:=\{(\sigma,n)\in 2^{<\omega}\times\omega \mid M\models K^{\mathbf{0}^{k}}_{R}(\sigma)\le n\}$ is still a member of the standard set in the proof of 1 $\to$ 2. No change is needed for the proof of 2 $\to$ 1.)
\end{proof}

\begin{theorem}\label{complex-in-arithmetic}
 Let $X\in 2^{\omega}$. Then, the following are equivalent.
\begin{enumerate}
 \item $X$ is $R$-complex relative to $\CPA$.
 \item For any recursive extension $T'\supseteq T^{R}$, $X$ is not compressible in $T'$, in other words, there exists $c\in\omega$ such that $T'\not\vdash K(X\rest n)<n-c$ for any $n\in\omega$.
 \item $X$ is not compressible in $T^{R*}$.
\end{enumerate}
\end{theorem}


On the other hand, randomness relative to CPA can be characterized by the notion of strong-randomness as follows.
For $A,B\subseteq 2^{<\omega}$, we write $A\prec B$ if for any $\sigma\in A$ there exists $\tau\in B$ such that $\tau\subseteq \sigma$.
\begin{defi}
 Let $m$ be a pre-measure.
 Then, we define a pre-measure $m^{*}$ as follows:
\begin{align*}
 m^{*}(F)&:=\inf\{m(C)\mid C\subseteq_{\mr{fin}} 2^{<\omega},~F\prec C\}.
\end{align*}
Then, $X\in 2^{\omega}$ is said to be \textit{strongly-$m$-random} if it is $m^{*}$-random.
\end{defi}
This definition agrees with the definition of strong-$h$-randomness ($\mr{pwt}_{h}$-randomness).
In fact, if $h$ is a convex recursive function, $(\mr{dwt}_{h})^{*}$-random is called vehement-$h$-random (named by Bj{\o}rn Kjos-Hanssen), and it is equivalent to the concept of strong-$h$-randomness (see Reimann~\cite{Reimann-APAL2008}).
More generally, for any recursive $h:2^{<\omega}\to\omega$, $(\mr{dwt}_{h})^{*}$-randomness is equivalent to $\mr{pwt}_{h}$-randomness by Theorem~\ref{gen-random-relPA} and the following theorem in \cite{propagation}: $X$ is $h$-random relative to CPA if and only if it is strongly-$h$-random.

\begin{lem}\label{Monotonicity-star}
Let $m$ be a pre-measure.
Then, for any $A,B\subseteq 2^{<\omega}$ such that $A\prec B$,
 we have $m^{*}(A)\le m^{*}(B)$.
\end{lem}
\begin{proof}
Trivial from the definition.
\end{proof}
The following theorem is a generalization of the $\CPA$ relativization theorem of Martin-L\"of randomness which was independently obtained by Downey, Hirschfeldt, Miller and Nies \cite[Proposition~7.4]{Downey-et-el-2005} and Reimann and Slaman \cite[Lemma~4.5]{Reimann-Slaman-2015}.
\begin{theorem}\label{gen-random-relPA}
Let $m$ be a pre-measure.
 Then, the following are equivalent.
\begin{enumerate}
 \item $X$ is strongly-${m}$-random.
 \item $X$ is ${m^{*}}$-random relative to $\CPA$. 
 \item $X$ is ${m}$-random relative to $\CPA$.
\end{enumerate}
\end{theorem}
\begin{proof}
We first show 1 $\to$ 2.
Assume that $X$ is not $m^{*}$-random relative to any CPA degree.
For each $Z\in2^{\omega}$, let $\{U^{Z}_{n}\mid n\in\omega\}$ be a universal $m^{*}$-test relative to $Z$.
Then, $X\in \bigcap_{Z\in \mr{CPA}}[U^{Z}_{n}]$ for any $n\in\omega$.
Define a uniformly r.e.~sequence $\{W_{n}\mid n\in\omega\}$ as $W_{n}=\{\sigma\in 2^{<\omega}\mid \A Z\in\CPA \E \tau\in U^{Z}_{n}(\tau\subseteq\sigma)\}$.
Then, $W_{n}\prec U^{Z}_{n}$ for any $Z\in\CPA$.
Since $\A Z\in \CPA\  \E k(X\rest k\in U^{Z}_{n})$, we have $\E k\  \A Z\in \CPA\  \E l<k(X\rest l\in U^{Z}_{n})$ by compactness of $\Pi^{0,X}_{1}$-class.
Thus, $X\in \bigcap_{n\in\omega}[W_{n}]$.
By Lemma~\ref{Monotonicity-star}, $\{W_{n}\mid n\in\omega\}$ is an $m^{*}$-test, thus $X$ is not $m^{*}$-random.

2 $\to$ 3 is trivial.

Finally we will prove 3 $\to$ 1.
 Assume that $X$ is $m$-random relative to $Z\in\mathrm{CPA}$ and is not $m^{*}$-random.
 Then, there exist a $m^{*}$-test $\{A_{n}\mid n\in\omega\}$ such that $X\in \bigcap_{i}[A_{i}]$.
 Now, we define a $\Pi^{0}_{1}$-class $P\subseteq 2^{<\omega}\times\omega$ such that
\begin{align*}
 W=\{ W_{n}\mid n\in\omega\}\in P \Leftrightarrow \A n (A_{n}\prec W_{n}\wedge m(W_{n})\le 2^{-n+1}).
\end{align*}
Note that $P$ is not empty since $\{A_{n}\mid n\in\omega\}$ is an $m^{*}$-test.
Hence, we can find $ W=\{ W_{n}\mid i\in\omega\}\in P$ such that $W\le_{T}Z$, and then $X\in \bigcap_{n\in\omega}[W_{n}]$.
Thus, $X$ is not $m$-random relative to $Z$, which is a contradiction.
\end{proof}

By Theorem~\ref{general-Schnorr}, we say that $X\in2^{\omega}$ is strongly-$R$-complex if it is $((R^{\surd})^{*})^{\surd}$-complex.
Now, we have the following by Theorems~\ref{complex-nonst}, \ref{complex-in-arithmetic} and \ref{gen-random-relPA}.
\begin{corollary}\label{cor-complex-in-arith}
Let $X\in 2^{\omega}$.
 Then, the following are equivalent.
\begin{enumerate}
 \item $X$ is strongly-$R$-complex.
 \item $X$ is $R$-complex relative to $\CPA$.
 \item $X$ is $R$-complex from a nonstandard model $M\models T$.
 \item For any recursive extension $T'\supseteq T^{R}$, $X$ is not compressible in $T'$.
\end{enumerate}
\end{corollary}

This also shows that a priori complexity $\mr{KA}$-$h$ works well even in an arithmetic, \textit{i.e.}, $X$ is $\mr{KA}$-$h$-complex if and only if it is not $\mr{KA}$-$h$-compressible in an arithmetic extending $T^{R_{\mr{KA}h}}$.
In particular, $X$ is Martin-L\"of random if and only if it is not $\mr{KA}$-compressible in an arithmetic extending $T^{R_{\mr{KA}}}$ if and only if it is not $\mr{KP}$-compressible in an arithmetic extending $T^{R_{\mr{KP}}}$, which is our claim in the first paragraph of this section.
(In that case, $R_{\mr{KP}}^{\surd}$ is just the usual fair-coin measure, and thus $R_{\mr{KP}}^{\surd}{}^{*}=R_{\mr{KP}}^{\surd}$ and strongly-$R_{\mr{KP}}$-complex is Martin-L\"of random, while Martin-L\"of random relative to $\CPA$ is again just Martin-L\"of random.)

On the other hand, $\mr{KP}$-$h$ does not work well in arithmetic in general.
By Reimann and Stephan\cite{Reimann/Stephan}, there exists a $1/2$-random $X$ which is not strongly-$1/2$-random.
(Here, (strongly-)$1/2$-random means (strongly-)$h$-random for $h(\sigma)=1/2|\sigma|$.)
Then, this $X$ is $\mr{KP}$-$1/2$-complex but it is $\mr{KP}$-$1/2$-compressible in some arithmetic extending $T^{R_{KP\text{-}1/2}}$.

Finally, we show that two different styles of definition of complexity coincide by relativization to ${\CPA}$.
\begin{prop}\label{R-and-hatR}
Let $R$ be a rule for a mode.
Then, $X\in2^{\omega}$ is $R$-complex relative to $\CPA$ if and only if it is $\hat{R}$-complex relative to $\CPA$.
\end{prop}
\begin{proof}
We show that if $X\in2^{\omega}$ is $R$-complex relative to $\CPA$, then it is $\hat{R}$-complex relative to $\CPA$.

Let $X\in 2^{\omega}$ be not $\hat{R}$-complex relative to $\CPA$, and let $Z\in \CPA$.
By Therems~\ref{general-Schnorr} and \ref{gen-random-relPA}, $X$ is not $(\hat{R}^{\surd})^{*}$-random.
Thus, there exists $(\hat{R}^{\surd})^{*}$-test $\{U_{i}\mid i\in\omega\}$ such that $X\in\bigcap_{i\in\omega}U_{i}$.
Define a $\Pi^{0}_{1}$-class $P\subseteq 2^{<\omega}\times 2^{<\omega}$ as
\begin{align*}
Y\in P \leftrightarrow& \E \{V_{i}\mid i\in\omega\} \A i (U_{i+1}\prec V_{i}\wedge R^{\surd}(V_{i})\le 2^{-i})\\
&\wedge \A s\in [2^{<\omega}\times\Z]^{<\omega}(s\subseteq \{(\sigma,|\sigma|-i)\mid \sigma\in V_{2i}\}\to \E r\subseteq_{\mr{fin}} Y(\ring{\hat{r}}=\ring{s}\wedge s\prec\hat{r}^{-2}\wedge r\in R)).
\end{align*}
By (bounded) K\"onig's lemma, this $P$ is non-empty.
Then, there exists $Y\in P$ such that $Y\le_{T}Z$.
Let $M^{Z}$ be an $R$-optimal mode relative to $Z$ and $K=K^{M^{Z}}$.
Then, there exists $c\in\omega$ such that for any $(\tau,\sigma)\in Y$, $K(\sigma)<|\tau|+c$.
By the definition of $P$, for any $i\in\omega$, there exist $\sigma\subseteq X$ and $\tau\in2^{<\omega}$ such that $(\tau,\sigma)\in Y$ and $|\tau|\le |\sigma|-i+2$, thus, $K(\sigma)\le |\sigma|-i+2+c$.
Hence, $X$ is not $R$-complex relative to $Z$.
\end{proof}
This shows that $\mr{KA}$-complex relative to $\CPA$  coincides with $\mr{KM}$-complex relative to $\CPA$.

\bibliographystyle{plain}
\bibliography{compbib}

\end{document}